\documentclass{article}

\usepackage{arxiv}
\usepackage{float} 
\usepackage[utf8]{inputenc} % allow utf-8 input
\usepackage[T1]{fontenc}    % use 8-bit T1 fonts
\usepackage{hyperref}       % hyperlinks
\usepackage{url}            % simple URL typesetting
\usepackage{booktabs}       % professional-quality tables
\usepackage{amsfonts}       % blackboard math symbols
\usepackage{nicefrac}       % compact symbols for 1/2, etc.
\usepackage{microtype}      % microtypography
\usepackage{lipsum}
\usepackage{graphicx}
\graphicspath{ {./images/} }
\usepackage{amsmath} 

\title{\textbf{How We Decide the Future of the Olympics ?}
}

\author{
 Wenlin Luo \\
  School of Chemistry and Chemical Engineering
  Nanjing University\\
  Nanjing 210023, P. R. China \\
  \texttt{luowenlin863@gmail.com} \\
   \And
 Chenghui Li \\
  School of Chemistry and Chemical Engineering
  Nanjing University\\
  Nanjing 210023, P. R. China \\
  \texttt{chli@nju.edu.cn}\\
}

\begin{document}
\maketitle
\begin{abstract}
The "Olympic Agenda 2020" stresses the urgency of measuring the impact of the Olympic Games on host cities in the face of declining bids. This paper presents the Hosting the Olympics Influence Evaluation Model (HOIEM) to assess these impacts and propose sustainable solutions.We select indicators (economic, socio-cultural, human, environmental and political) based on literature review and construct HOIEM. We use AHP and TOPSIS-EWM to determine the final weights for these indicators. We identify 45 potential host cities based on IOC requirements. For the Winter Olympics, secondary screening and the GM(1,1) model highlight Calgary, Canada as the top city. For the Summer Olympics, the SWOT analysis identifies Beijing, China. We propose to hold Spring, Summer, Autumn and Winter Olympics every 4 years, with fixed cities for Summer and Winter and bidding for Spring and Autumn.Finally, we conduct sensitivity analysis for HOIEM, demonstrating the independence and relevance of selected indicators through semi-quantitative visual analysis. 
\end{abstract}

% keywords can be removed

\section{Introduction}
As the number of countries bidding for various Olympics sharply declines, the IOC gradually realizes the modern ethical dilemmas that the Olympics face today\cite{1}. Looking back at the hosting processes of the Olympics in recent decades, the host countries that successfully held the Games not only did not obtain significant macroeconomic benefits as expected, but also incurred actual bidding costs that were several times higher than the declared costs and "Olympic legacy" that could not be handled. These facts indicate that significant reforms for the Olympics are urgently needed, as reflected in the "Olympic Agenda 2020+5"[1]. Therefore, establishing a reasonable evaluation model to guide the exploration of various new strategies and policies is crucial for improving the attractiveness and sustainability of the Olympics.

\section{Hosting the Olympics Influence Evaluation Model (HOIEM) }
\label{sec:headings}

\subsection{Indicator selection}
\begin{figure}[H]
    \centering
    \includegraphics[width=0.5\linewidth]{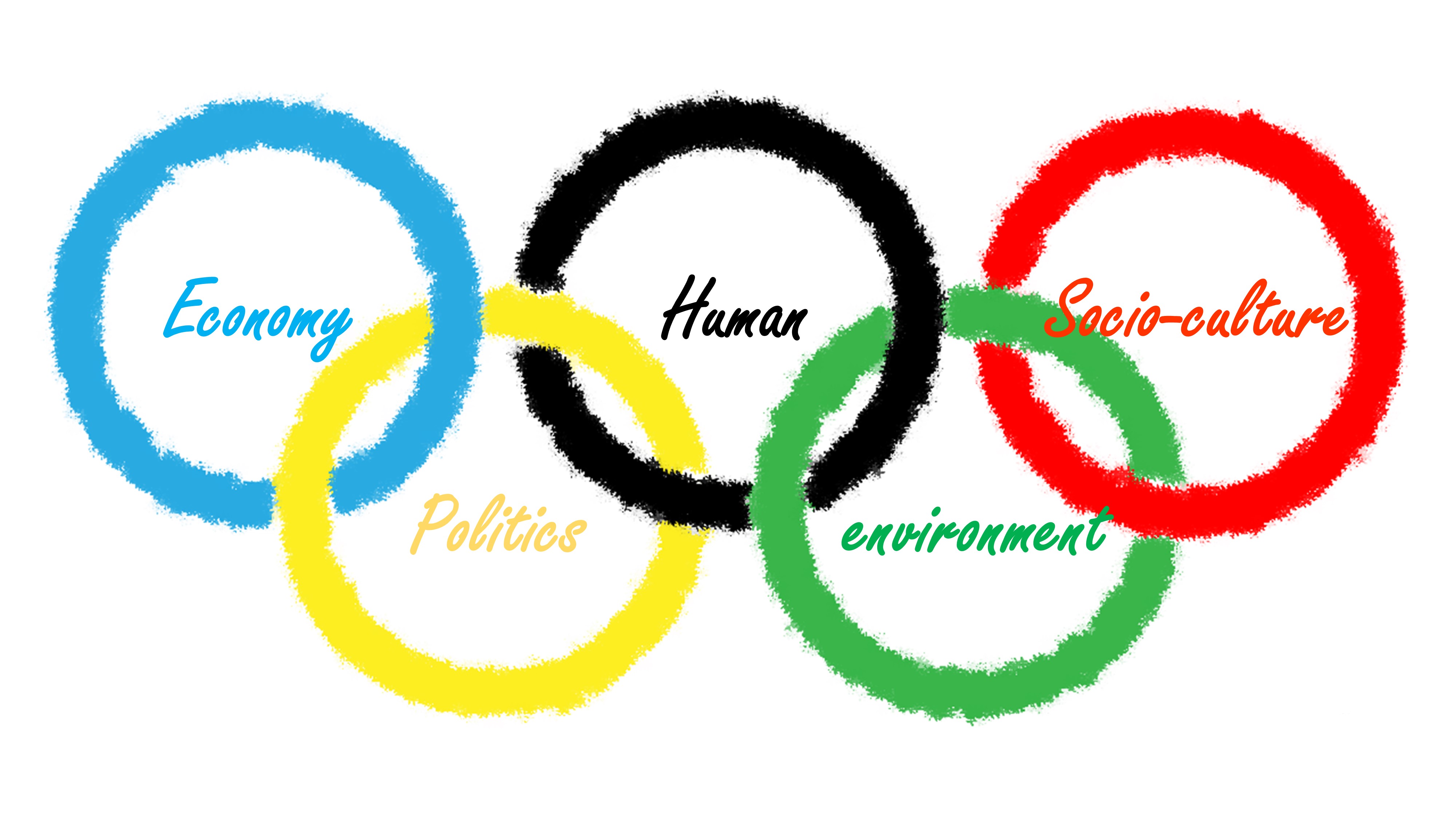}
    \caption{Evaluation Factors}
    \label{1}
\end{figure}

There are many factors that contribute to the impact of the Olympics on host cities. When selecting factors to establish an evaluation model for the impact of hosting the Olympics, more evaluation indicators are not necessarily better. Based on the original question and references, we have summarized \textbf{5 }aspects of impact evaluation factors, including Economy, Human, Socio-culture, Politics, and Environment\textbf{ }(\textbf{as shown in Figure}\ref{1} ), to help us build a comprehensive evaluation system and find corresponding data support.
\   \newline
From the numerous literature on the analysis of the impacts of the Olympics listed in \textbf{Table 1}, we extracted the following secondary indicators (where positive and negative factors are identified by the symbols (+) and (-), respectively). These indicators were selected and categorized based on the characteristics of data composition and the requirements of a comprehensive, universal, and minimally overlapping evaluation model.
\newline
Based on the review of the literature listed in\textbf{ Table} \textbf{2}, we developed a multidimensional evaluation model for the impacts of hosting the Olympics. The model can comprehensively reflect the various impacts of hosting the sporting event on the host city, providing quantitative basis for addressing the challenges faced by the Olympics. According to this model, we selected specific secondary indicators guided by the 5 primary indicators, in order to use concise and non-overlapping indicators to maximize the expression of the meaning of each primary indicator.
\newline
\begin{table}[]\arraystretch{1.5}
\caption{Extraction of secondary indicators}
\label{c1}

\begin{tabular}{lll}
\hline
\textbf{Primary indicators} & \textbf{Secondary indicators}                                                                                                                                                                     & \textbf{Researchers}                                                                       \\ \hline
Economy& \begin{tabular}[c]{@{}l@{}}Employment rate\\ GDP per capita\\ Income from events\end{tabular}                                                                                                     & \begin{tabular}[c]{@{}l@{}}
Arthur, \cite{2} 
Andresen, M.A., \& Tong, W \cite{3}
\end{tabular}\\ \hline
Human                       & \begin{tabular}[c]{@{}l@{}}Athlete evaluation\\ Audience evaluation\\ Other people evaluation\end{tabular}                                                                                        & \begin{tabular}[c]{@{}l@{}}
Lee Ludvigsen\cite{4} 
Grix, J., \& Lee, D\cite{5}
\end{tabular}\\ \hline
Sociocultural               & \begin{tabular}[c]{@{}l@{}}Strengthen project execution ability\\ Carry forward the Olympic spirit\\ Enhancement of international image\\ Widening gap between the rich and the poor\end{tabular} & \begin{tabular}[c]{@{}l@{}}
Jarvis, J\cite{6} 
Horne, J., \& Manzenreiter, W\cite{7}
\end{tabular}\\ \hline
Political                   & \begin{tabular}[c]{@{}l@{}}For the host region\\ For international\end{tabular}                                                                                                                   & \begin{tabular}[c]{@{}l@{}}
Staiano\cite{8}
\end{tabular}\\ \hline
Environmental               & \begin{tabular}[c]{@{}l@{}}Promote environmental protection\\ Generate a large amount of waste\end{tabular}                                                                                       & \begin{tabular}[c]{@{}l@{}}
Romanazzi, V., \& Rovere, R\cite{9}
\end{tabular}\\ \hline
\end{tabular}
\end{table}
\subsubsection{Economic Indicators}

\begin{figure}[H]
    \centering
    \includegraphics[width=0.8\linewidth]{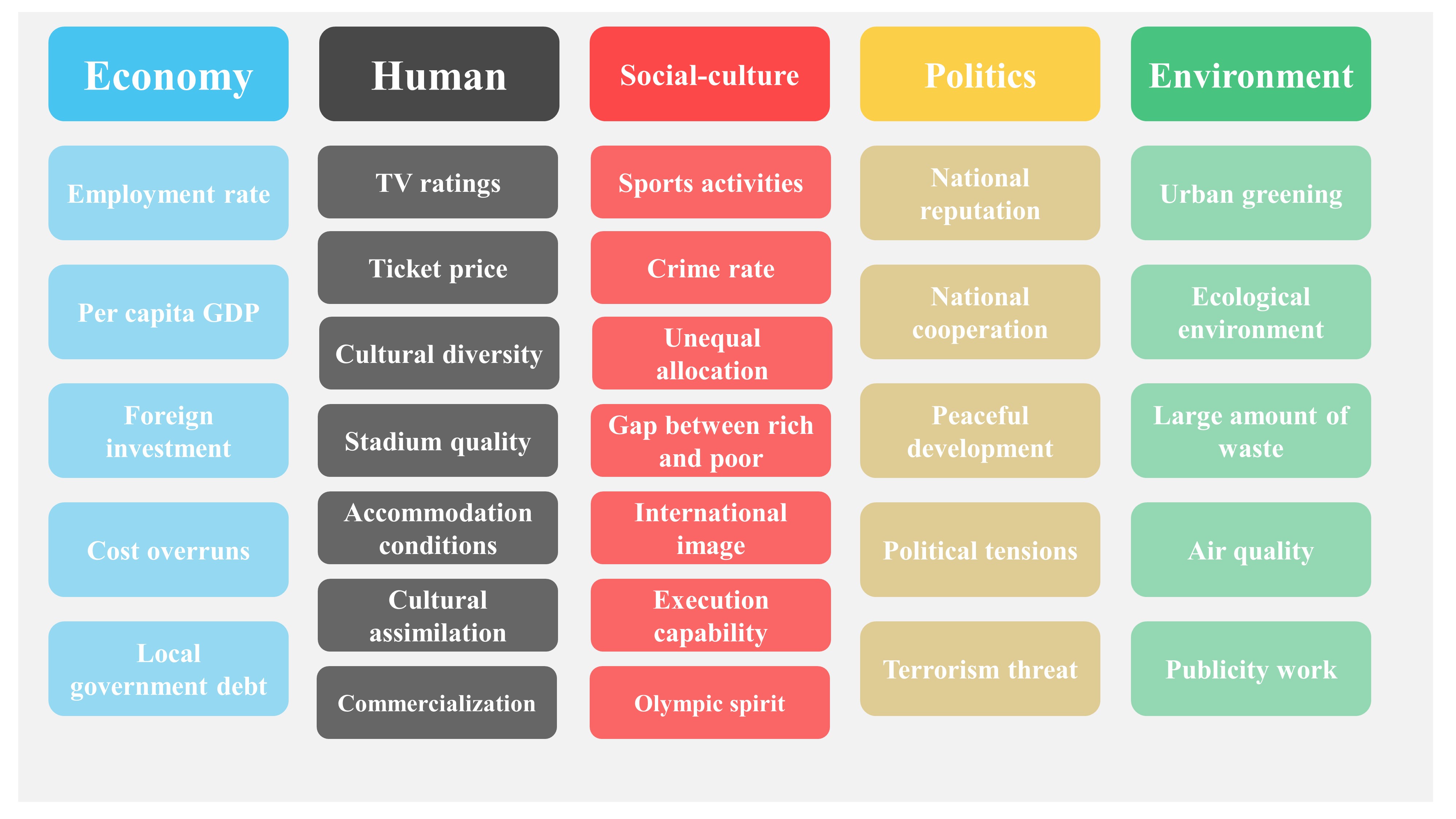}
    \caption{Secondary Indicators}
    \label{2}
\end{figure}

As shown in \textbf{Figure} \ref{2}, we can see that in terms of economic factors, based on the secondary indicators obtained from reviewing the literature, we have explored the data composition characteristics of economic types and divided them into six secondary indicators (which we can refer to as \textit{A1 }to\textit{ A6}). It can be determined that the secondary indicator system we constructed can comprehensively and systematically include the various economic impacts of hosting the Olympics on the host country. For example, not only does it discuss the macroeconomic growth brought by hosting the Olympics to the host city, but also considers the huge economic costs of hosting the Olympics. In addition, while maintaining breadth, we maximize the use of data characteristics to maintain the objectivity of the secondary mechanism of action on the primary indicators.

\subsubsection{Human Indicators}

Similar to the process of determining the secondary indicators for economic factors, after determining the interaction mechanism for human factors' secondary indicators, we also conducted a secondary analysis and ultimately divided them into seven secondary indicators (denoted as \textit{B1 }to \textit{B7}).

\subsubsection{Sociocultural,Political and Environmental Indicators}

To maintain the generality of the model, we have noticed that only discussing the impact of the Olympics on the economy and humans is far from sufficient. For example, the political terrorist incident at the 1972 Munich Olympics and the impact of COVID-19 on the 2020 Tokyo Olympics are both black swan events, but they can have a huge impact on the host. At the same time, the promoting role of the Olympic spirit in social culture seems to be significant. Finally, we continue to pay attention to reports on irreversible damage to the ecological environment and the living environment of local residents before and after the Olympics. Based on this, we hope to extract some secondary indicators that can reflect these aspects, which is conducive to establishing a more scientific evaluation system. They are respectively denoted as \textit{C1 }to \textit{C7},\textit{D1} to \textit{D5} and \textit{E1} to\textit{ E5}.

\subsection{Establishment and solution of HOIEM}
\subsubsection{Determination of Weight}
The Analytic Hierarchy Process (AHP), developed by Saaty in 1980\cite{10}, is a comprehensive evaluation method that combines systems analysis and decision-making. It helps us to take into account both qualitative and quantitative factors when making decisions and then arrive at a comprehensive result.Firstly, the relationships between various factors are analyzed and a systematic hierarchical structure is established.Secondly, pairwise comparisons of the elements within the same level are made, and judgment matrices are constructed.
\begin{equation}b_{ii}=1,b_{ij}=b_{ji}^{-1}>0\end{equation}
In equation (1), $b_{ii}$ and $b_{ij}$ represent different evaluation indicators.
Based on the judgment matrix, the relative weights of the compared elements to the criterion are calculated, and consistency check is performed.
\begin{equation}CI=\frac{\lambda_{\max}-n}{n-1}\end{equation}

\begin{equation}CR=\frac{CI}{RI}\end{equation}
The equation above shows that $ \lambda_{max}$ is the maximum eigenvalue of the n-order judgment matrix,  $CI$ is the consistency index, and $RI$ is the corresponding average random consistency index. If the calculated consistency ratio $CR<0.1$, the consistency check of the judgment matrix is considered to pass, otherwise it needs to be revised. 
Finally, the corresponding weight values $V_{j}$ are calculated based on the judgment matrix and solved using the analytic hierarchy process.
\newline
TOPSIS (C.L. Hwang and K. Yoon, 1981)\cite{11}, which stands for Technique for Order Preference by Similarity to an Ideal Solution, is a comprehensive evaluation method that fully reflects the original information. Modified by EWM, it achieves objective weighting based on the variability of each indicator. First, each data in the original matrix is normalized in a positive direction.
\begin{equation}M=\max\left\{a-\min\left\{x_{i}\right\},\max\left\{x_{i}\right\}-b\right\}\end{equation}
\begin{equation}
\tilde{x}_i =
\begin{cases}
1 - \frac{a - x_i}{M}, & x_i < a \\
1, & a \leq x_i \leq b \\
1 - \frac{x_i - b}{M}, & x_i > b
\end{cases}
\end{equation}
\begin{equation}z_{ij}=\frac{x_{ij}}{\sqrt{\sum_{i=1}^nx_{ij}^2}}\end{equation}
Based on the compositional characteristics of the collected data, we select equation (6) as the standardization method. Then, we solve the probabilities used in the relative entropy calculation\cite{12}.
\begin{equation}p_{ij}=\frac{\widetilde{z}_{ij}}{\sum_{{i\operatorname{=}1}}^{n}\widetilde{z}_{ij}}\end{equation}
In formula (7), $p_{ij}$ represents the relative entropy probability of the $j$-th sample under the $i$ -th index. Finally, the information entropy of the index is calculated to determine the entropy weight.
\begin{equation}e_j=-\frac1{\ln n}\sum_{i\operatorname{=}1}^np_{ij}\ln\left(p_{ij}\right)\end{equation}
\begin{equation}H_j=\frac{1-e_j}{n-\sum_je_j}\end{equation}
The equation above shows that $e_{j}$ represents the information entropy contained in the $j$-th indicator, and $H_{j}$ represents the entropy weight carried by the $j$-th indicator.
\newline
The combination weighting model based on objective modification of subjectivity can well compensate for the "superposition effect" and "multiplication effect" caused by traditional additive and multiplicative combination weighting methods. The combined weights can reflect the expert's intention to a certain extent while embodying the information carried by the sample data itself. First, calculate the importance ratio between adjacent indicators. 
\begin{equation}\breve{s}_{j}=\sqrt{\frac{1}{m}\sum_{{i\operatorname{=}1}}^{m}\left(x_{ij}-\overline{{x_{ij}}}\frac{H_{j}+V_{j}}{H_{j}V_{j}}\right)^{2}}\end{equation}
\begin{equation}
r_k = \left\{
\begin{array}{c}
\min\left\{2,\frac{\breve{s}_{k-1}}{\breve{s}_k}\right\}, \breve{s}_{k-1} \geq \breve{s}_k \\
1, \breve{s}_{k-1} < \breve{s}_k
\end{array}
\right.
\end{equation}
In formula (11),  $s_{j}$ represents the normalized standard deviation of the attribute $j$ , and its importance ratio is denoted as $r_{k}$ .
\begin{equation}W_m=\left(1+\sum_{k=2}^m\prod_{j=k}^mr_j\right)^{-1}\\W_{j-1}=r_jW_j\end{equation}
Finally, the comprehensive weight of the indicator to the criterion layer can be determined according to the above formula.

\subsubsection{Solution of HOIEM}
\begin{figure}[H]
    \centering
    \includegraphics[width=0.8\linewidth]{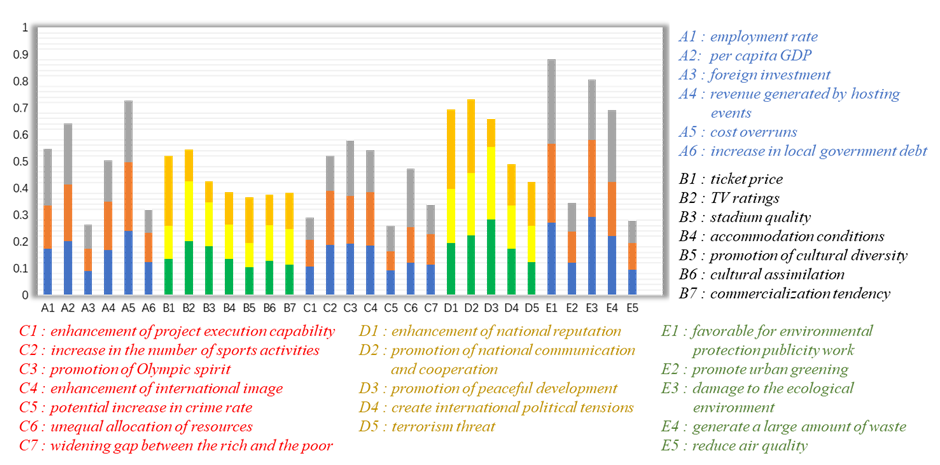}
    \caption{Comprehensive weight of secondary indicators}
    \label{3}
\end{figure}

Based on the previous work, we obtained the subjective and objective weights of the corresponding indicators through the AHP and EWM-TOPSIS methods, and obtained the comprehensive weights through the model method based on objective correction of subjectivity. The weight values of different indicators calculated are shown in \textbf{Figure \ref{3}}.

we determined the weights u for the five primary indicators and independent subjective, objective, and comprehensive weights for each secondary indicator under each primary indicator. The total weight is calculated as follows.
\begin{equation}\Omega_{ij}=U_iu_j\end{equation}
Calculation of the weights of each secondary index is shown in\textbf{ Table 3}.
We used the idea in principal component analysis that the sum of the top 10 comprehensive weights was 0.735. Therefore, in evaluating the generality and simplicity of the model, we selected these 10 characteristic factors  $\xi_{j}$ for further analysis. We found that they came from different categories of primary indicators, which is also a manifestation of the preservation of generality.After secondary weighting normalization of the selected characteristic factors, we obtained the adjusted weights  $\gamma_{j}$ , as shown in Table \ref{b2}.

\begin{table}[]
\renewcommand\arraystretch{1.5}
\label{b2}
\caption{Feature indicators and its weight}
\begin{tabular}{lllllllllll}
\hline
Characteristic factor number & $\xi_{1}$& $\xi_{2}$& $\xi_{3}$& $\xi_{4}$& $\xi_{5}$& $\xi_{6}$& $\xi_{7}$& $\xi_{8}$& $\xi_{9}$& $\xi_{10}$\\ \hline
$\gamma_{j}$& 0.165 & 0.146 & 0.132 & 0.124 & 0.105 & 0.088 & 0.082 & 0.071 & 0.059 & 0.030 \\ \hline
Original number              & A5    & A4    & C1& C7    & A2    & E4    & B2    & D2    & C3    & D5    \\ \hline
\end{tabular}
\end{table}

Finally, we obtain the evaluation formula for the impact assessment model of hosting the Olympic Games as shown in equation (14).
\begin{equation}\chi=\sum_{j=1}^{10}\gamma_j\xi_j\end{equation}
Formula (14) defines the evaluation function ${\chi}$ for the impact assessment of hosting the Olympics. It can be seen that ${\chi}$ is a positive function defined in the background of considering important factors, including:
\newline
\textbf{·The feature indicator group has generality and comes from different five primary indicators.}
\newline
\textbf{·The feature indicators maintain independent relationships with each other without overlapping.}
\newline
\textbf{·The selection of feature indicators combines subjective and objective weights with high credibility.}
\newline
\textbf{·The function describes well that the cost issue in the economic aspect is the most important factor that the Olympics brings impacts on.}
\newline
\textbf{·It also takes into account unpredictable factors that may affect the hosting of the Olympics.}
\newline
\textbf{·The evaluation function confirms that the economic and sociocultural aspects are the main contributions to the impact assessment model of the Olympics.}

\section{Application of the HOIEM}
\label{sec:others}
\begin{figure}[H]
    \centering
    \includegraphics[width=0.8\linewidth]{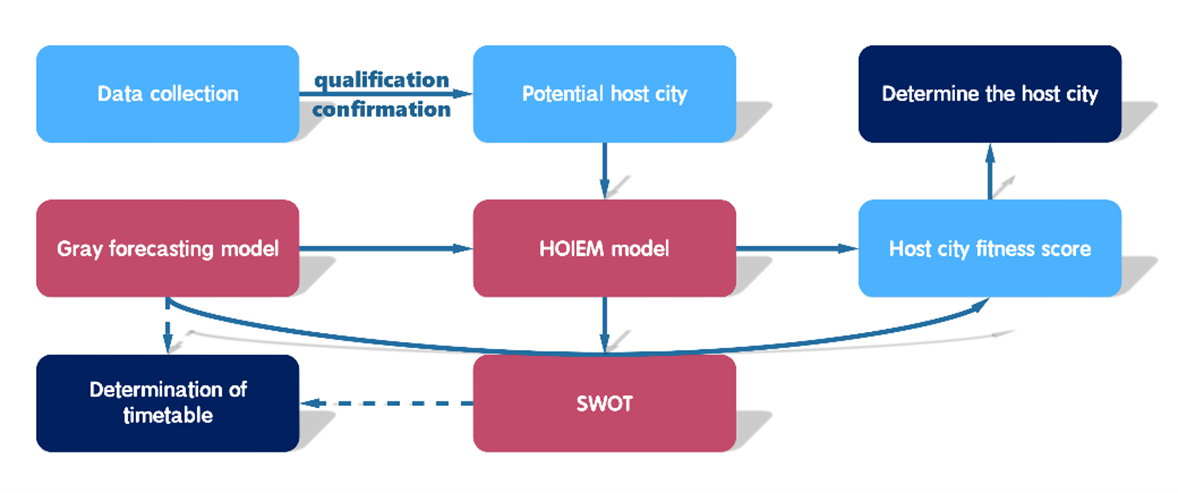}
    \caption{Olympics Comprehensive Decision Model}
    \label{4}
\end{figure}

In order to propose innovative decisions to address the difficulties faced by cities in bidding for the Olympic Games.Consider using the GM (1,1) model, combined with the HOIEM model constructed in the first question, to obtain the potential changes in characteristic factors of candidate cities before and after hosting the Olympics, thereby determining their adaptability scores for hosting the Olympics. As shown in \textbf{Figure \ref{4}}, based on the gray prediction system, we need to combine the characteristics of the secondary indicators extracted for different types in the HOIEM model to obtain the corresponding Olympic host city fitness score, which determines the potential for a city to host the Olympics. Finally, guided by the HOIEM model, we use the GM (1,1) model to predict the hosting schedule.
\subsection{Screening of Potential Host Cities}
According to the standardization recommendations of the International Olympic Committee, we used the total GDP of a country and its level of sports development as the first screening criteria, and \textbf{selected 45 potential candidate cities} from the countries that ranked high in the comprehensive results. The results are shown in \textbf{Table 3}. Based on the HOIEM model, we believe that only potential candidate cities meet the most basic conditions for hosting the Olympic Games, so the actual host city must be a non-empty proper subset of the potential candidate city group.

\begin{table}[]\centering
\renewcommand\arraystretch{1.5}
\label{b3}
\caption{45 potential host cities}
\begin{tabular}{lllll}
\hline
New York                  & Shenzhen   & Lyon     & Madrid        & Brussels     \\ \hline
Tokyo   Metropolitan Area & Seoul      & Turin    & Mexico   City & Bangkok      \\ \hline
Los   Angeles             & Moscow     & Rome     & Jakarta       & Dublin       \\ \hline
Shanghai                  & Osaka      & Toronto  & Amsterdam     & Jerusalem    \\ \hline
London                    & Mumbai     & Calgary  & Riyadh        & Buenos Aires \\ \hline
Paris                     & Delhi      & Busan    & Ankara        & Oslo         \\ \hline
Beijing                   & Birmingham & Moscow   & Bern          & Vienna       \\ \hline
Chicago                   & Berlin     & Brasília & Warsaw        & Abuja        \\ \hline
Philadelphia              & Frankfurt  & Sydney   & Stockholm     & Copenhagen   \\ \hline
\end{tabular}
\end{table}
\subsection{Determine the host city of the Winter Olympics}
GM(1,1) is a single-variable grey prediction model with a first-order difference equation\cite{13}. It has advantages such as requiring a small sample size and having a simple model establishment and solution process. The specific implementation process is as follows:Given a time series $X^{(0)}=\{X^{(0)}(1),X^{(0)}(2),...,X^{(0)}(n)\}$ with $n$ observations, a new sequence $X^{(1)}=\{X^{(1)}(1),X^{(1)}(2),...,X^{(1)}(n)\}$ is generated by accumulating the sum.The corresponding differential equation of GM(1,1) can be expressed as the following formula:
\begin{equation}\frac{dX^{(1)}}{dt}+\alpha X^{(1)}=\mu \end{equation}
In the GM(1,1) model, $\alpha $ is called the developing grey degree and $\mu$ is called the internal control grey degree. The general solution to the differential equation is:
\begin{equation}\widehat{X}^{(1)}(k+1)=\left[X^{(0)}(1)-\frac\mu\alpha\right]e^{-\alpha k}+\frac\mu\alpha \end{equation}
Based on the preliminary screening of potential host cities, we further identify the most suitable city to host the Winter Olympics. Before predicting factors such as average temperature and snowfall for candidate host cities, we first use common knowledge from geography and climate to eliminate 36 cities that are obviously unsuitable for hosting the Winter Olympics. The various characteristics listed in the table all indicate that they are not suitable to be included in the list of candidate host cities for the Winter Olympics. After determining 9 potential candidate host cities, we used the GM(1,1) model to pre-dict their future (until 2050) average temperature and snowfall in February to determine if they meet the International Olympic Committee's requirements for hosting the Winter Olympics, which include an average temperature below 0℃ in February, an ideal temperature range of -17℃ to -10℃, and a minimum snowfall of 30cm in February. After this work was completed, we eliminated two-thirds of the potential host cities. Finally, the remaining 3 cities seem to have no difference in terms of the rigid requirements. In order to better distinguish them, we can only use the hosting suitability score function constructed earlier under the guidance of the HOIEM model.
\newline
We will input the predicted results into the HOIEM evaluation model and solve the corresponding scoring function. The results are as follows: Moscow S\textsubscript{1}=1.381(S\textsubscript{base}=0.6,S\textsubscript{evaluate}=0.781), Pyeongchang S\textsubscript{2}=1.602(S\textsubscript{base}=0.8,S\textsubscript{evaluate}=0.802), and Calgary S\textsubscript{3}=1.649(S\textsubscript{base}=0.8,S\textsubscript{evaluate}=0.849). It can be seen that Calgary, Canada is the most ideal location for hosting the Winter Olympics.

\subsection{Determine the host city of the Summer Olympics}
\subsubsection{Secondary screening of potential host cities}
\begin{figure}[H]
    \centering
    \includegraphics[width=0.6\linewidth]{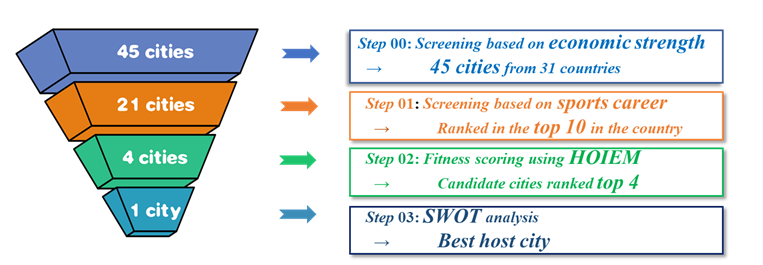}
    \caption{Evaluation and selection process of host cities for the Summer Olympics}
    \label{5}
\end{figure}

For the Summer Olympics, we can naturally follow the layer-by-layer progressive screening method for the host cities of the Winter Olympics. The difference is that the International Olympic Committee's requirements for the host city of the Summer Olympics do not seem to have quite rigid conditions like the Winter Olympics. For this reason, we slightly changed the screening process this time, as shown in \textbf{Figure \ref{5}}.

The measurement standard of sports career in the above picture comes from the medal points ranking of the previous Summer Olympics from 1896 to 2020 (including 5 points for gold medals, 1 point for silver medals, and 0.5 points for bronze medals). The process of solving the fitness score of the host city in HOIEM is similar to that of the Winter Olympics and will not be repeated here.

\subsubsection{SWOT analysis}
\begin{figure}[H]
    \centering
    \includegraphics[width=1\linewidth]{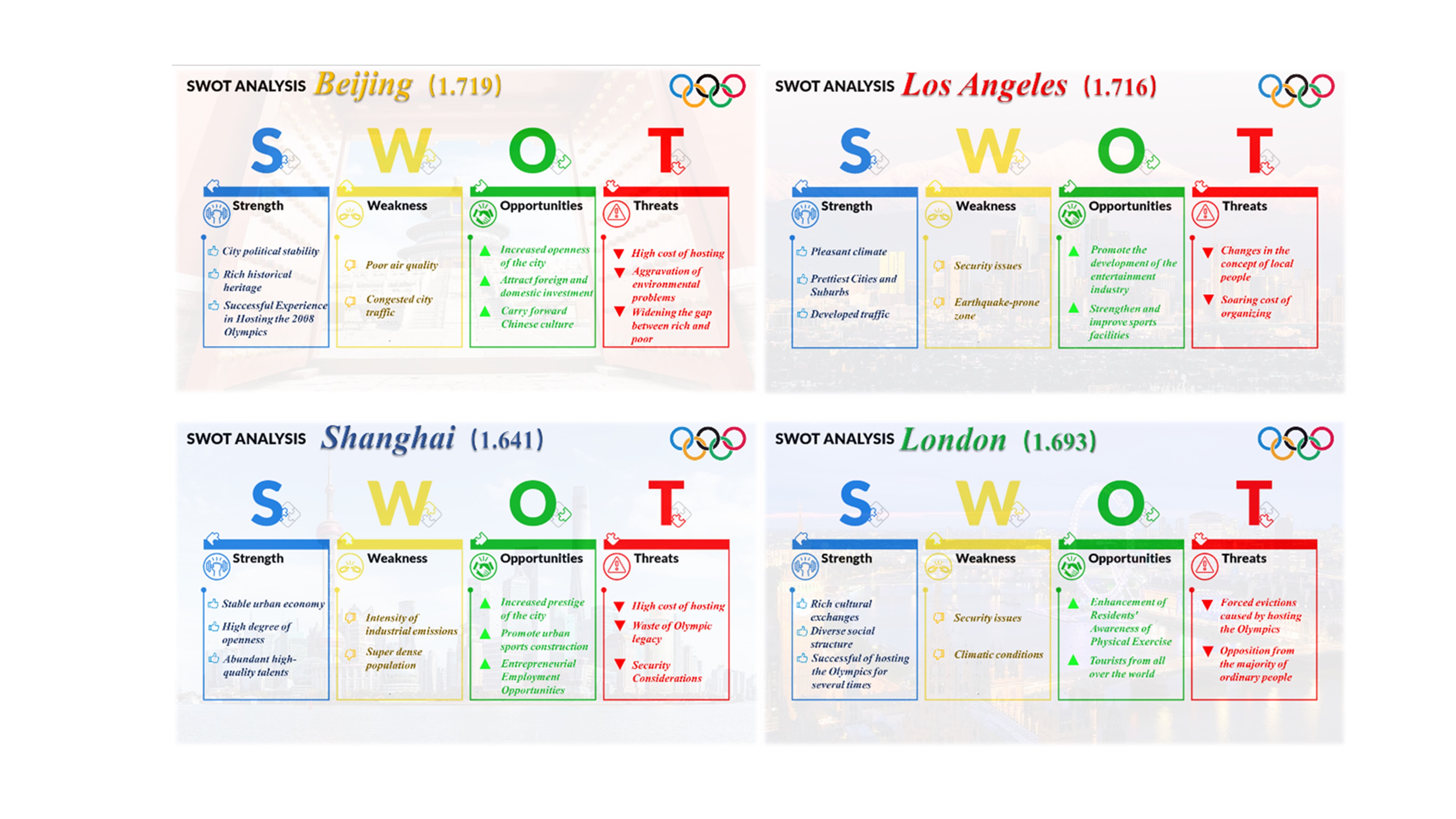}
    \caption{Analysis of 4 Alternative Cities}
    \label{6}
\end{figure}

Finally, we choose the SWOT analysis method. It is an analytical method used to identify the strengths, weaknesses, opportunities and threats of the decision-making object. We performed SWOT analysis on the four cities obtained earlier, and the results are shown in the \textbf{Figure \ref{6}}.

Based on the SWOT analysis and the fitness score results in HOIEM, we believe that Beijing is the most suitable city for hosting the Summer Olympics.\textbf{ }As mentioned in the conclusion of the case study on the Beijing 2008 Olympic Games in (4.4.2) above. We have reasons to believe that choosing Beijing as the fixed venue for the Olympic Games can bring a grand sports event to people all over the world while controlling the cost of hosting the Olympic Games as much as possible.
In this section of the analysis on the Winter and Summer Olympic Games, we have obtained the most suitable cities for hosting them, which are Calgary and Beijing respectively. Next, we discuss how the Olympic Games are held, while quantitatively analyzing its merits using our established HOIEM.

\subsection{Determine the manner of holding the future Olympics}
On the basis of combining the topic description and consulting relevant literature, the following four innovative solutions are perfectly summarized.Plan A believes that on the basis of not changing the holding time, the Summer and Winter Olympics should have corresponding permanent venues. Plan B indicates that 4 Olympics in spring, summer, autumn and winter should be held in each cycle (four years), and they should all have corresponding fixed venues. Plan C does not support having a fixed venue, and only proposes that four Olympic Games should be held in each cycle. Plan D is a com-promise between Plans B and C. It will be held four times in a cycle and the Summer and Winter Olympics will have fixed venues.
Next, we discuss the above four possible schemes quantitatively and semi-quantitatively based on the model technology established above. And compare it with the original way of hosting the Olympics, draw their respective advantages and disad-vantages, and determine the best plan as much as possible.

\subsubsection{Quantitative analysis based on HOIEM}
\begin{figure}[H]
    \centering
    \includegraphics[width=0.5\linewidth]{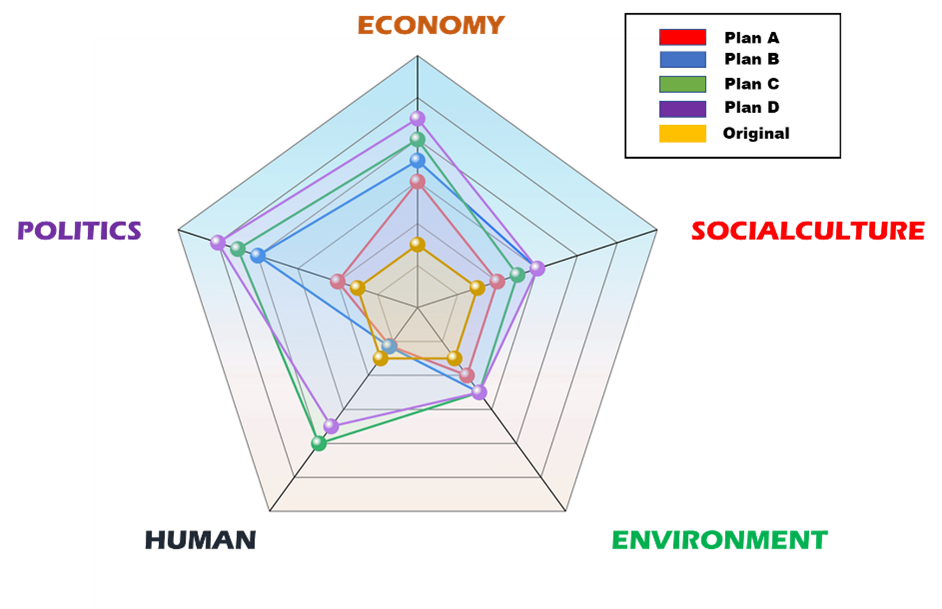}
    \caption{Evaluation results of 4 schemes on characteristic groups}
    \label{7}
\end{figure}

In order to intuitively and clearly show the effects of the four schemes on the city compared with the original holding method, we consider using the characteristic index group of HOIEM for analysis. Evaluate its impact on each feature index separately, and finally get the total impact function and compare its size. To this end, we use the comparison method used in the establishment of the AHP model to define an objective function for evaluating the impact, as shown in the \textbf{Table 4}:

\begin{table}[]\centering
\renewcommand\arraystretch{1.5}
\label{b4}
\caption{Objective Evaluation Scale of Influence }
\begin{tabular}{llllll}
\hline
Impact Scale & 1         & 3                  & 5               & 7                   & 9                    \\ \hline
illustrate   & no effect & slightly favorable & quite favorable & extremely favorable & absolutely favorable \\ \hline
\end{tabular}
\end{table}
With this definition, we can evaluate the improvements of the four schemes for the ten characteristic group marks used in HOIEM combined with the original scheme. Obviously, the greater the degree, the better. The evaluation results and their visualizations are shown in  \textbf{Figure \ref{7}}.
we can see that Scheme D seems to win in terms of integration.\textbf{ }It is worth mentioning that the economic, ecological, environmental and political aspects of the Plan A, \textit{B}, \textit{C}, and \textit{D} all have positive changes compared with the original scheme, which shows that the direction of this improvement is in line with science. Perhaps on the basis of Plan D, we can also consider a better plan. However, due to space limitations, we will not make further program setting and evaluation. At the end, we will conduct a SWOT analysis of the corresponding scheme in order to comprehensively analyze the characteristics of all aspects of the scheme.

\section{Sensitivity Analysis of HOIEM}
\subsection{Sensitivity Analysis of obsolete indicators}

\begin{figure}[H]
    \centering
    \includegraphics[width=0.6\linewidth]{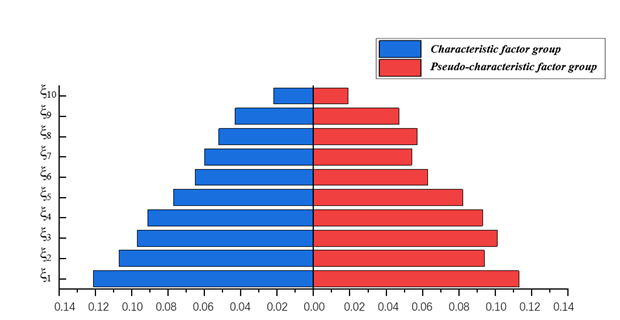}
    \caption{Sensitivity analysis results of the impact assessment model}
    \label{8}
\end{figure}

In the evaluation model, sensitivity analysis is usually used to analyze the fluctuation of the selection of different characteristic factors for the evaluation model. In the previous work, we identified 10 characteristic factors, and in the process of sensitivity analysis, we considered randomly selecting 5 out of 20 unselected factors among the 30 candidate factors as False feature factors, the random 5 of the feature factors are replaced equally, the result is shown in \textbf{Figure \ref{8}}.
\subsection{Semi-quantitative-Visual Analysis based on RSM}
Box-Behnken Response Surface Method (BBD-RSM)\cite{14} is a commonly used tool in optimization visualization analysis, and can also be used to evaluate the stability analysis of the model. The basic principle is to use the weight of a certain factor as the system response value and the total evaluation function as the global variable, and use an intuitive graphic language to display its interaction relationship, so as to judge the rationality of the selected factors under semi-quantitative conditions.
\begin{figure}[H]
    \centering
    \includegraphics[width=0.8\linewidth]{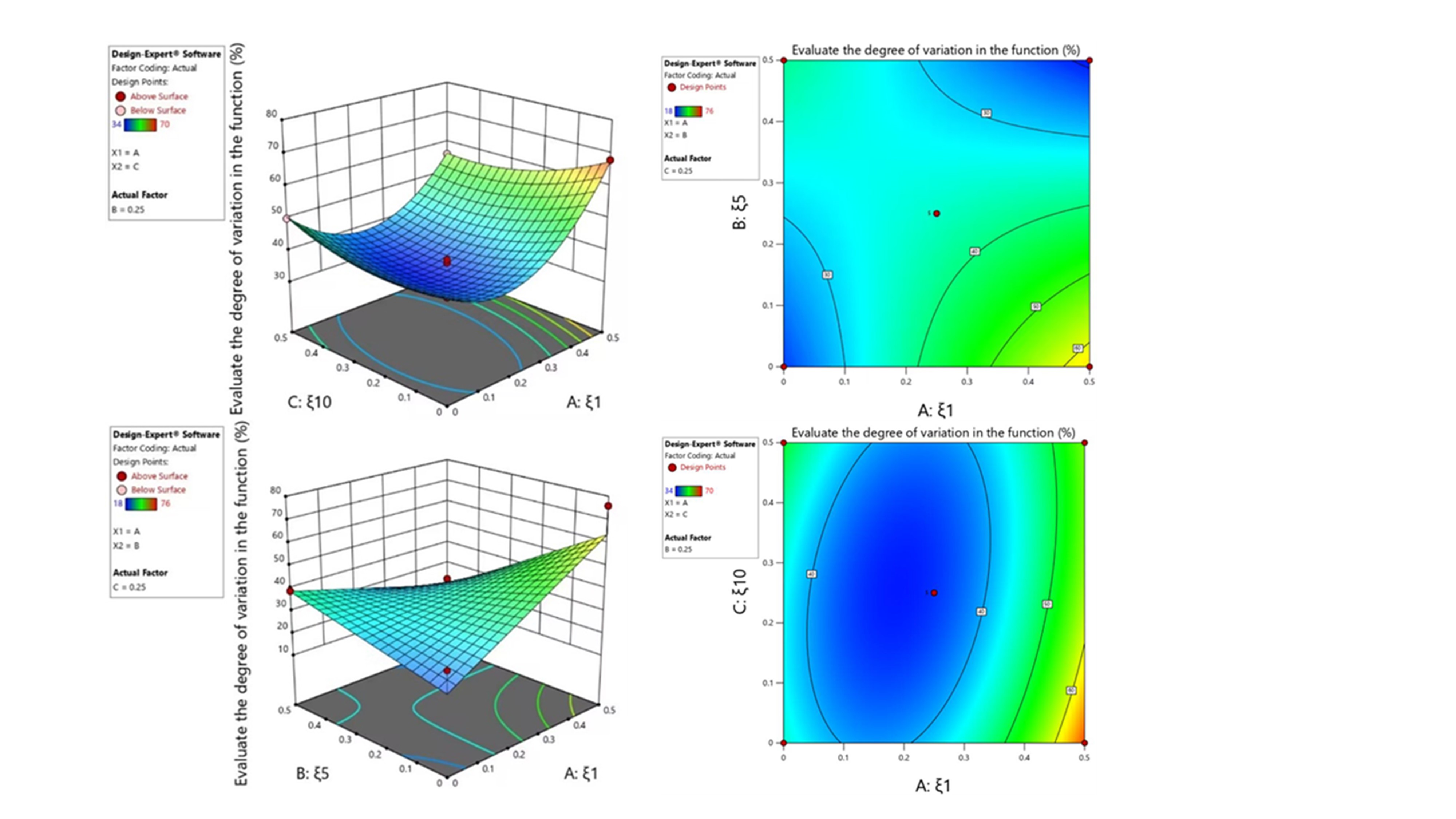}
    \caption{Response surfaces of the evaluation function with respect to $\xi_{1}$ and $\xi_{10}$ }
    \label{9}
\end{figure}

Analysis of \textbf{Figure} \ref{9} shows that the weight changes of  $\xi_{1}$ and $\xi_{10}$ will cause obvious response to the variation degree of the evaluation function. It can be seen that the response extreme value brought by $\xi_{1}$ is about 69\%, and the influence extreme value of $\xi_{10} $ is about 50\%. This result supports the representative of feature index selection. At the same time, the response extreme value brought by $\xi_{1}$ and $\xi_{10}$ together is about 25\%, which is lower than the two independent response extreme values and is the global minimum, indicating that these two characteristic factors are independent\textbf{.}Similar to the above analysis, as shown in Figure 20, we also tested the responsiveness of  $\xi_{1}$ and $ \xi_{5}$to the degree of variation of the evaluation function.

\section{Conclusion}
We adopted a combined weighting model based on AHP-EWM, incorporating both objective correction and subjective weight determination, and selected more comprehensive indicators to establish the HOIEM (Host City Olympic Impact Evaluation Model). In the sensitivity analysis of HOIEM, we utilized a combination of quantitative and semi-quantitative-visualization methods; the two methods complement each other, fully illustrating the stability and feasibility of the HOIEM model. To ensure that the article is quite readable, we have adopted a large number of visualization methods for the solution results and analysis of the model. Our HOIEM is based on five key aspects: economy, socio-cultural, people, politics, and ecological environment, comprehensively including all aspects of the influence of the Olympic Games on the host city, thus demonstrating its good practicality and promotion value.
\begin{figure}[H]
    \centering
    \includegraphics[width=0.8\linewidth]{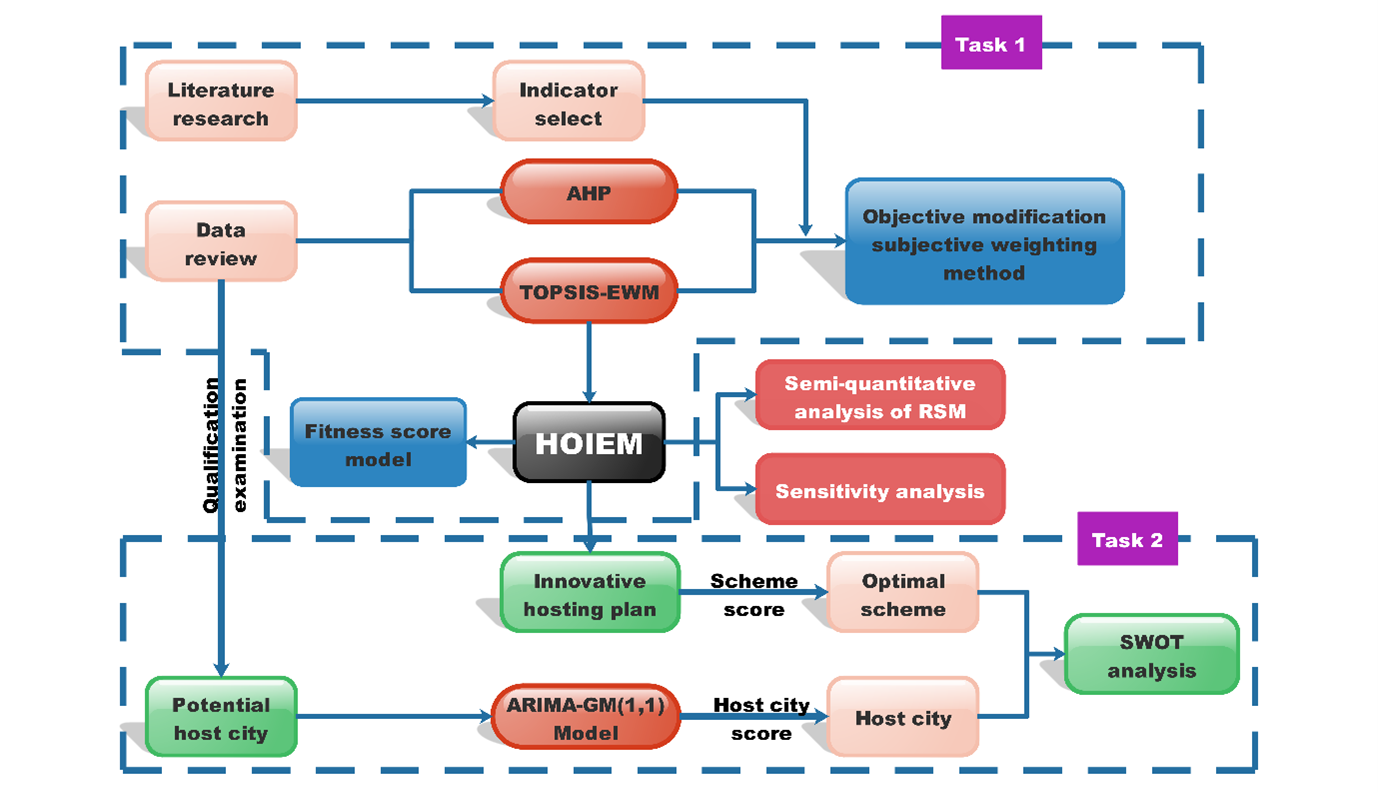}
    \caption{Our work}
    \label{10}
\end{figure}

\bibliographystyle{IEEEtran}
\small\bibliography{references}

\end{document}